\documentclass[12pt]{amsart}
\usepackage{amscd,amssymb}
\theoremstyle{latex 2e}
\newtheorem{thm}[subsection]{Theorem}
\newtheorem{lem}[subsection]{Lemma}

\newtheorem{con}[subsection]{Conjecture}
\newtheorem{rem}[subsection]{Remark}
\newtheorem{defn}[subsection]{Definition}

\numberwithin{equation}{section}

\newcommand{\CC}{{\mathbb C}}
\newcommand{\R}{{\mathbb R}}

\newcommand{\M}{\mathfrak{M}}

\begin{document}

\title[K\"ahler manifolds with nef Ricci class] %
{K\"ahler manifolds with numerically effective Ricci class and
maximal first Betti number are tori}

\author{Fuquan Fang}
\address{Nankai Institute of Mathematics,
Weijin Road 94, Tianjin 300071, P.R.China}
\email{ffang@nankai.edu.cn}

\thanks{**Supported by NSFC Grant 19925104 and 973 Project for Foundation Science
of China}

\begin{abstract}
Let $M$ be an $n$-dimensional K\"ahler manifold with numerically
effective Ricci class. In this note we prove that, if the first
Betti number $b_1(M)=2n$, then $M$ is biholomorphic to the complex
torus $T^n_\CC$.
\end{abstract}
\maketitle

\section{Introduction}
\label{sec:intro}

\baselineskip=.25in

\vskip2mm

Let $M$ be a compact complex manifold with a fixed hermitian
metric $\omega$. By [DPS93] [DPS94] a holomorphic line bundle $L$ over $M$
is called {\it numerically effective} (abb. nef) if for every
$\varepsilon >0$, there is a smooth hermitian metric
$h_\varepsilon$ on $L$ such that the curvature satisfies
$$\Theta _{h_\varepsilon}\ge - \varepsilon\omega$$
If $M$ is projective, $L$ is nef precisely if $L\cdot C\ge 0$ for
all curves $C\subset M$. We say a K\"ahler manifold $M$ is nef if
the anticanonical bundle $-K_M$ is nef. In [DPS93] it is
conjectured, for a nef K\"ahler manifold, both of the following
holds:

(A1) the fundamental group $\pi _1(M)$ has polynomial growth.

(A2) the Albanense map $\alpha : M\to Alb(M)$ is surjective.

If $M$ is projective, (A2) was proved by Q. Zhang [Zh]. In [Pa1]
Paun proved (A1), assuming a conjecture of Gromov concerning the
fundamental group of Riemannian manifold with almost non-negative
Ricci curvature (compare [ChCo]).

By  the Aubin-Calabi-Yau theorem [Ab][Ya], [DPS93] proved that $M$
is nef if and only if there exist a sequence of K\"ahler metrics
$\{\omega _k\}$ on $M$ such that, for each $k>0$, the metric
$\omega _k$ belongs to a fixed cohomology class $\{\omega\}$, and
the Ricci curvature of $\omega _k$ is bounded from below by
$-1/k$.

A Bochner type theorem for the first Betti number was obtained by
Paun [Pa2], namely, for every nef K\"ahler manifold $M$ of complex
dimension $n$ it holds that $b_1(M)\le 2n$.  The main result of
this note is the following:

\vskip 2mm

\begin{thm} Let $M$ be a nef K\"ahler manifold of dimension $n$. If the
first Betti number $b_1(M)=2n$, then $M$ is biholomorphic to the
complex torus $T^n_\CC$.
\end{thm}

\vskip 2mm

Theorem 1.1 may be considered as a complex version of a conjecture
of Gromov, proved by T. Colding [Co], asserts that a Riemannian
$n$-manifold of almost non-negative Ricci curvature and first
Betti number $n$ is homeomorphic to the torus $T^n$.

Obviously, Theorem 1.1 implies conjecture (A2) in the case of
$b_1(M)=2n$. In the proof of Theorem 1.1, in fact we will prove
that there is a uniform upper bound for the diameters. But this
does not hold if $b_1(M)=2n-2$ (the first Betti number of a
K\"ahler manifold is always even). Since there is a sequence of
K\"ahler metrics on $S^2$ in the same K\"ahler class with positive
Ricci curvature but converge to a non-compact space of dimension
$1$, thus the product $T^{n-1}_\CC\times S^2$ serves as an
example.

The following result verifies (A2) for manifold with
$b_1(M)=2n-2$, provided $G_1/[G_1, G]$ has rank at least two,
where $G=\pi _1(M)$, $G_1=[G, G]$.

\vskip 2mm

\begin{thm} Let $M$ be a nef K\"ahler manifold of dimension $n$.
Let $G=\pi _1(M)$. If the first Betti number $b_1(M)=2n-2$, and
$G_1/[G_1, G]$ has rank at least two where $G_1=[G, G]$, then the
Albanese map $\alpha : M\to T^{n-1}_\CC$ is surjective.
\end{thm}

\vskip 2mm

\begin{rem}

By Theorems 1.1 one confirms immediately conjecture (A2) for
$n=2$. This was first obtained in [DPS93] using algebraic geometry
methods.
\end{rem}

\vskip 2mm

The proof of our Theorems uses the deep results in Riemannian
geometry, including the equivariant Gromov-Hausdorff convergence
[FY], a splitting theorem of Cheeger-Colding for limit spaces
[ChCo], and a result of Paun [Pa2]. It would be interesting if
Theorem 1.1 could be proved in a way of pure algebraic geometry.
Indeed, if the Albanese map $\alpha$ is surjective, by the
Poincare-Lelong equation, one obtains easily that  a nef K\"ahler
manifold $M$ of dimension $n$ with $b_1(M)=2n$ is biholomorphic to
the complex torus $T^n_\CC$ (compare [Mo]).

By Campana [Ca] the above conjectures (A1) and (A2) together with
Gromov's celebrated theorem [Gr] implies that the fundamental
group of a nef K\"ahler manifold is almost abelian. By our
approach, it seems plausible to prove the following:

\vskip 2mm

\begin{con} Let $M$ be a nef K\"ahler manifold of dimension $n$. If 
there is an epimorphism $\varphi: \pi _1(M)\to \Gamma$ where $\Gamma$ is a
torsion free nilpotent group of rank at least $2n$, then $\Gamma
\cong \Bbb Z^{2n}$ and $M$ is biholomorphic to the complex torus
$T^n_\CC$.
\end{con}

\vskip4mm

{\bf Acknowledgement:} The author would like to thank Ngaiming Mok
for helpful discussions concerning the Poincare-Lelong equation. The paper is finished during 
author's visit to the Max-Plack Institut f\"ur Mathematik. The author is very grateful
to the institute for financial support.

\vskip 8mm

\section{Proof of Theorems 1.1 and 1.2}

\vskip  8mm

By [DPS93], a nef K\"ahler manifold $M$ admits a family of
K\"alher metrics $\omega _\varepsilon$ in the same K\"ahler class
$[\omega ]$ with Ricci curvature $\text{Ric.} (\omega _\varepsilon
)\ge -\varepsilon \omega $, where $\varepsilon \in (0, 1)$. The
diameters of this family may not have a uniform upper bound. In
other words, the pointed Gromov-Hausdorff limit of $(M, \omega
_\varepsilon)$ may not be compact. Because of this, many
techniques in metric geometry do not apply to this situation. To
overcome this difficulty, [DPS93] obtained the following key
lemma.

\begin{lem}[DPS93]
Let $M$ be a nef K\"ahler manifold. Let $U\subset \tilde M$ (the
universal covering of $M$) be a  compact subset. Then, $\forall
\delta >0$, there exists a closed subset $U_{\varepsilon,
\delta}\subset U$ such that

(2.1.1) $\text{vol} _\omega (U- U_{\varepsilon, \delta})<\delta$;

(2.1.2) $\text{diam}_{\omega _\varepsilon}(U_{\varepsilon,
\delta})\le C/\sqrt \delta$.

\noindent where $C$ is a  constant independent of $\varepsilon$
and $\delta$.

\end{lem}

\vskip 2mm

For convenience let us recall the definition of equivariant
Gromov-Hausdorff distance (cf. [FY] for details).

Let $\M$ (resp.  $\M _{eq}$) denote the set of all isometry
classes of pointed metric spaces $(X, p)$ (resp. triples $(X,
\Gamma, p)$), such that, for any $D$, the metric ball $B_p(D, X)$
of radius $D$ is relatively compact and such that $X$ is a length
space [GLP][FY] (resp.  $(X, p)\in \M$ and $\Gamma$ is a closed
subgroup of isometries of $X$).

Let $\Gamma (D)=\{\gamma \in \Gamma : d(\gamma p, p)<D\}$.

\begin{defn} Let $(X, \Gamma , p), (Y, \Lambda , q)\in \M _{eq}$.
An $\varepsilon$-equivariant pointed Hausdorff approximation
stands for a triple $(f, \varphi, \psi)$ of maps $f: B_p(\frac 1
\varepsilon, X)\to Y$, $\varphi : \Gamma (\frac 1 \varepsilon)\to
\Lambda (\frac 1 \varepsilon)$, and $\psi :\Lambda (\frac 1
\varepsilon)\to \Gamma (\frac 1 \varepsilon)$ such that

(2.2.1) $f(p)=q$;

(2.2.2) the $\varepsilon$-neighborhood of $f(B_p(\frac 1
\varepsilon, X))$ contains $B_q(\frac 1 \varepsilon,Y)$;

(2.2.3) if $x, y\in B_p(\frac 1 \varepsilon , X)$, then
$$
|d(f(x), f(y))-d(x, y)|<\varepsilon;
$$

(2.2.4) if $\gamma \in \Gamma (\frac 1 \varepsilon)$, $x\in
B_p(\frac 1 \varepsilon, X)$, $\gamma x\in B_p(\frac 1
\varepsilon, X)$, then
$$
d(f(\gamma x), \varphi (\gamma )(f(x)))<\varphi ;
$$

(2.2.5) if $\mu \in \Lambda (\frac 1 \varepsilon)$, $x\in
B_p(\frac 1 \varepsilon, X)$, $\psi (\mu )(x)\in B_p(\frac 1
\varepsilon, X)$, then
$$
d(f(\psi (\mu)(x)), \mu f(x))<\varepsilon.
$$
\end{defn}

The {\it equivariant pointed Gromov-Hausdorff distance}
$d_{eGH}((X,\Gamma, p), (Y, \Lambda , q))$ stands for the infimum
of the positive numbers $\varepsilon$ such that there exist
$\varepsilon$-equivariant pointed Hausdorff approximations from
$(X, \Gamma, p)$ to $(Y, \Lambda , q)$ and from  $(Y, \Lambda ,
q))$ to $(X, \Gamma, p)$.

\begin{proof}[Proof of Theorem 1.1] Let $\omega _k$ be a sequence
of K\"ahler metrics on $M$ in the same K\"ahler class with Ricci
curvature $\ge -\frac 1k \omega $. Let $\tilde M_k$ be the
Riemannian covering space of $M_k$ (the manifold $M$ with the
K\"ahler metric $\omega _k$). Using Lemma 2.1 Paun [Pa2] proved
that there is an open subset $\tilde U_k\subset \tilde M_k$ of
$\mbox{diam}_{\omega _k}(\tilde U_k )\le C$ such that the
homomorphism $\pi _1(U_k)\to \pi _1(M_k)$ is surjective, where
$U_k$ is the image of $\tilde U_k$ in $M_k$, $C$ is a universal
constant independent of $k$.

For convenience let $G=\pi _1(M)$, and let $\Gamma =G/[G, G]$.
Consider $\bar M_k=\tilde M_k/[G, G]$. By assumption $\Bbb
Z^{2n}\subset \Gamma $ acts on $\bar M_k$ by isometry. By a lemma
of Gromov [GLP] (compare [Pa 2]) it follows that there is a finite
index torsion free subgroup $\Gamma _k$ of $\Gamma$ such that,
fixing a base point $p_k \in \bar U_k\subset \bar M_k$,

(2.3.1) the geometric norm of any non-trivial element of $\Gamma_k
$ is at least $C$.

(2.3.2) $\Gamma_k $ is generated by $\gamma _1, \cdots ,\gamma
_{2n}$ so that the geometric norm of every $\gamma _i$ is at most
$2C$.

Since $\Gamma _k$ acts on $\bar M _k$ by isometry, the quotient
space $\bar M _k/\Gamma _k$ is a finite Riemannian covering space
of $M_k$. Because the Ricci curvature of $\bar M_k$ is bounded
from below, by the Gromov compactness theorem (cf. [FY]) the
pointed spaces converge
$$(\bar M_k, \Gamma _k,
p_k)\stackrel{d_{eGH}}\longrightarrow (X, \Gamma _\infty, q)$$ in
the equivariant Gromov-Hausdorff topology when $k$ tends to
infinity. By (2.3.1) it is easy to see that the isometric action
of $\Gamma _\infty$ on $X$ is discrete and effective.  By the
splitting theorem [ChCo] the limit space $X=Y\times \R^\ell$,
where $Y$ contains no line. By [GLP] it is well-known the
Hausdorff dimension of $X$ is at most $2n$, therefore $\ell\le
2n$. We first need

\vskip 2mm

\begin{lem}
$\Gamma _\infty\cong \Bbb Z^{2n}$.
\end{lem}

\begin{proof}[Proof of Lemma 2.3]

By definition, there are maps $\varphi _k: \Gamma _k(k)\to \Gamma
_\infty (k)$, $\psi _k: \Gamma _\infty (k) \to \Gamma _k(k)$ and a
$\frac 1k$-Hausdorff approximation $f_k: B_{p_k}(k, \bar M_k)\to
B_q (k, X)$ satisfying (2.2.1)-(2.2.5).

We first claim that $\varphi _k$ is injective for sufficiently
large $k$. If not, there are two elements $\gamma _k \ne \lambda
_k \in \Gamma _k(k)$ such that $\varphi (\gamma _k)=\varphi
(\lambda_k)$ for any $k$. Let $\mu _k=\varphi (\gamma _k)=\varphi
(\lambda_k)$. Put $x=p_k$. By (2.2.4) we get that $d(f_k(\lambda
_kx), \mu _kf_k(x))<\frac 1k$ and $d(f_k(\gamma _kx), \mu
_kf_k(x))<\frac 1k$. Therefore, $d(f_k(\lambda _kx), f_k(\gamma
_kx))<\frac 2k$ and so $d(\lambda _k\gamma _k^{-1}x, x)=d(\lambda
_kx, \gamma _kx)<\frac 4k$ since $f_k$ is a $\frac 1k$-Hausdorff
approximation. A contradiction to (2.3.1).

Secondly, we claim that $\varphi _k(\gamma _i\gamma _j)=\varphi
_k(\gamma _i)\varphi _k(\gamma _j)=\varphi _k(\gamma _j)\varphi
_k(\gamma _i)$ for any $\gamma _i, \gamma _j\in \Gamma _k(k)$ so
that $\gamma _i\gamma _j\in \Gamma _k(k)$. In fact, by (2.2.4)
again we get that $d(\varphi _k(\gamma _i\gamma
_j)f_k(x),f_k(\gamma _i\gamma _jx))<\frac 1k$; $d(\varphi
_k(\gamma _i)\varphi _k(\gamma _j)f_k(x), \varphi _k(\gamma
_i)f_k(\gamma _jx))<\frac 1k$ and $d(f_k(\gamma _i\gamma _jx),
\varphi _k(\gamma _i)f_k(\gamma _jx))<\frac 1k$. Thus, $d(\varphi
_k(\gamma _i\gamma _j)f_k(x), \varphi _k(\gamma _i)\varphi
_k(\gamma _j)f_k(x))<\frac 3k$. For the same reason as above, by
(2.3.1) it follows that $\varphi _k(\gamma _i\gamma _j)=\varphi
_k(\gamma _i)\varphi _k(\gamma _j)$. The claim follows.

Similar argument applies to show that $\varphi _k(\gamma
_i^{-1})=\varphi _k(\gamma _i)^{-1}$, if $\gamma _i,\gamma
_i^{-1}\in \Gamma _k(k)$.

Next we verify that $\varphi _k: \Gamma _k(k)\to \Gamma _\infty
(k)$ is also surjective.

We argue by contradiction. Assume such an element $\mu _k\in
\Gamma _\infty (k)$. By (2.2.5) $d(f_k(\psi (\mu _k)(x), \mu
_kf_k(x))<\frac 1k$. By (2.2.4) $d(f_k(\psi (\mu _k)(x),\varphi
_k(\psi( \mu _k))f_k(x))<\frac 1k$. Therefore, $d(\varphi _k(\psi(
\mu _k))f_k(x), \mu _kf_k(x))<\frac 2k$. By (2.3.1) this implies
that $\mu _k=\varphi _k(\psi( \mu _k))$. A contradiction.

For sufficiently large $k$, let $\Gamma _0$ be the subgroup of
$\Gamma _\infty$ generated by $\varphi_k(\gamma _1), \cdots ,
\varphi_k(\gamma _{2n})$. It may be verified easily that this does
not depend on the choice of $k$. By (2.3.2) and the above $\Gamma
_0$ is a commutative group of rank $2n$. Since $\varphi_k$ is
surjective, $\Gamma _0=\Gamma _\infty$. The desired result
follows.
\end{proof}

To continue the proof of Theorem 1.1, we first prove that
$X=\R^{2n}$. It suffices to show that $\ell =2n$.

We argue by contradiction. Assume $\ell<2n$.

Since $\Gamma _\infty$ preserves the splitting, there is a
well-defined homomorphism $p: \Gamma _\infty \to \mbox{Isom}(\R
^\ell)$. Let $\Gamma _{0,\infty}$ denote the kernel of $p$. By the
generalized Bieberbach  theorem (cf. [FY]) the image $p(\Gamma
_\infty)$ has rank at most $\ell$. By Lemma 2.3 $\Gamma
_{0,\infty}$ has rank $\ge 1$. For a non-trivial element of $\mu
\in \Gamma _{0,\infty}$, by (2.2.5) there is a sequence of element
$\gamma _k =\psi _k(\mu )\in \Gamma _k$ (of infinite order) such
that the $\gamma _k$-action on $\bar M_k$ converges to the action
of $\mu$ on $Y\times \Bbb R^\ell$. Observe that a minimal geodesic
representation in $\bar M_k/\Gamma_k$ gives rise a line in $\bar
M_k$, on which $\gamma _k$ acts by deck transformation. This
sequence of lines converges to a line in $Y$ on which $\mu$ acts
by translation. Therefore the line lies in $Y$ since $\mu \in
\Gamma _{0,\infty}$ acts trivially on the factor $\R ^\ell$. A
contradiction to the assumption that $Y$ has no line. Hence $\ell
={2n}$.

Finally, by (2.3.2) we see that $\R^{2n}/\Gamma _\infty$ is
compact. By [FY] Lemma 3.4 $\bar M_k/\Gamma _k$ converges to
$\R^{2n}/\Gamma _\infty$. This shows that $\bar M_k/\Gamma _k$ has
uniformly bounded diameter. Therefore, $\bar M_k/\Gamma _k$ has
almost non-negative Ricci curvature in Gromov's sense [GLP]. By
[Co] we conclude that $\bar M_k/\Gamma _k$ is homeomorphic to a
torus $T^{2n}$, and so is $M$. By Poincare-Lelong equation it
follows that the Albanese map has no zeros and is actually a
biholomorphism. This completes the proof of Theorem 1.1.
\end{proof}

\begin{rem}
The above proof actually shows that a sequence of K\"ahler metrics
on $T^{2n}$ in the same K\"ahler class $[\omega]$ with Ricci
curvature $\ge -\varepsilon \omega $ has uniformly bounded
diameter, and so the metrics do not collapse.
\end{rem}

\vskip  8mm

Let $G=\pi _1(M)$. Consider the lower central series 
$$\cdots G_2\subset G_1\subset G_0=G$$
where $G_1=[G, G]$ and $G_2=[G_1, G]$. 
Let $G_2'\subset G$ be the normal subgroup such that $G/G_2'=(G/G_2)/\text{torsion}$.
Assume $H_1(G)/\text{torsion}\cong \Bbb Z^{2n-2}$, and $\text{rank}(G/G_2')=2n-2+m$.
By [Pa2] we may assume elements 
$\gamma _1, \cdots, \gamma _{2n-2}; \alpha _1, \cdots , \alpha _m \in G$ which generate 
a finite index subgroup  $\Gamma _k'\subset G/G_2'$ and satisfy (2.3.1), (2.3.2) and

(2.5.1) the geometric norms of $\alpha _1, \cdots ,\alpha _m$ are
all less than $2C$.

We warn that it is not true if we require that $\alpha _1, \cdots ,\alpha _m$ 
satisfy (2.3.2).

Now we start the proof of Theorem 1.2. We will only sketch the main steps since the proof
follows the same line as the previous one.

\begin{proof}[Proof of Theorem 1.2]

Let $\bar M_k'=\tilde M_k/G_2'$. Consider the triple $(\bar M_k', \Gamma _k', p_k)$.
The pointed spaces converge 
$$(\bar M_k', \Gamma _k', p_k)\stackrel{d_{eGH}}\longrightarrow
(X, \Gamma _\infty ',q)$$

Exactly the same argument in the previous proof implies that $X=Y\times
\Bbb R^{2n-2}$ and $Y$ contains at least a line since the group generated by
$\{\alpha_1, \cdots ,\alpha _m\}$ converges to a non-trivial isometry group acting 
on $X$ acting trivially on $\Bbb R^{2n-2}$, where $Y$ is a length space of Hausdorff dimension
at most two. However, since (2.3.2) is not satisfied for the $\alpha _i$'s, the limit group
$\Gamma _\infty'$ may not be discrete (compare [FY] example 3.11). Therefore, $X=Y_0\times \Bbb R^{2n-1}$ where the
Hausdorff dimension of $Y_0$ is at most $1$.

If $Y_0$ is compact, e.g., zero dimensional, by
(2.3.2) and (2.5.1) it follows that the limit space $Y_0\times \Bbb R^{2n-1}/\Gamma _\infty'$ is 
compact. Therefore, the diameters of the sequence $\bar M_k'/\Gamma _k'$
have a uniform upper bound, so are the diameters of $M_k$ (since $M_k$ is a finite isometric 
quotient of $\bar M_k'/\Gamma _k'$).
By [Pa3] it follows that the Albanese map is surjective.

If $Y_0$ is $1$-dimensional and non-compact, clearly, $Y_0$ has two ends and thus $Y_0$ contains
a line. By [ChCo] once again $Y_0=\Bbb R$. This proves that $X\cong \Bbb R^{2n}$. 
Since $m\ge 2$, the rank of $\Gamma _\infty'$ is at least $2n$ (may be non-discrete). By the 
generalized Bieberbach theorem the quotient $\Bbb R^{2n}/\Gamma _\infty'$ has to cocompact.
For the same reasoning as above the desired result follows.
\end{proof}

\vskip 8mm

\centerline{\large References}

\vskip 8mm

[Ab] T. Aubin, {\em \'Equantions du type Monge-Amp\`ere sur les
vari\'et\'es k|'ahl\'erienne compactes}, Bull. Sci. Math. France,
{\bf 102} (1978), 63-95

[Ca] F. Campana, {\em Remarques sur les groupes de K\"ahler nilpotents,} {\bf 28}(1995), 
307-316

[ChCo] J. Cheeger; T. Colding, {\em Lower bounds on Rcci curvature
and almost rigidity of wraped products}, Ann. Math.,  {\bf 144}
(1996), 189-237

[Co] T. Colding, {\em Ricci curvature and volume convergence},
Ann. Math., {\bf 145}(1997), 477-501

[DPS 93] J.P. Demailly, T. Peternell, M. Schneider, {\em K\'ahler
manifolds with numerically effective Ricci class}, Comp. Math.,
{\bf 89} (1993), 217-240

[DPS 94] J.P. Demailly, T. Peternell, M. Schneider, {\em Compact
complex manifolds with numerically effective tangent bundles}, J.
Alg. Geometry, {\bf 3} (1994), 295-345

[FY] K. Fukaya; T. Yamaguchi, {\em The fundamental groups of
almost non-negatively curved manifolds}, Ann. of Math. {\bf 136}
(1992),  253-333

[Gr] M. Gromov, {\em Group of polynomial growth and expanding
maps}, Publ. Math. I.H.E.S., {\bf 53}(1981), 53-73

[GLP] M. Gromov, J. Lafontaine; P.Pansu, {\em Structures metriques
pour les varietes riemannienes},  CedicFernand Paris, 1981

[Mo] N. Mok, {\em Bounds on the dimension of $L^2$-holomorphic sections of vector 
bundles over complete Kähler manifolds of finite volume}, Math. Z. {\bf 191}(1986), 303--317

[Pa1] M. Paun, {\em Sur le groupe fondamental des vari\'et\'es
k\"ahl\'eriennes compactes \`a classe de Ricci num\'eriquement
effective}, C. R. Acad. Paris, {\bf t. 324} (1997), 1249-1254

[Pa2] M. Paun, {\em Sur vari\'et\'es k\"ahl\'eriennes compactes
\`a classe de Ricci num\'eriquement effective}, Bull. Sci. Math.,
{\bf 122} (1998), 83-92

[Pa3] M. Paun, {\em On the Albanese map of compact K\"ahler
manifolds with numerically effective Ricci curvature}, Comm. Anal.
Geom. {\bf 9}(2001), 35-60

[Ya] S. T. Yau, {\em On the Ricci curvature of a complex K\"ahler
manifold and the complex Monge-Amp\'ere equation I}, Comm. Pure
and Appl. Math., {\bf 31}(1978), 339-411

[Zh] Q. Zhang, {\em On projective manifolds with nef anticanonical
bundle}, J. reine. angew. Math., {\bf 478}(1996), 57-60

\end{document}